\documentclass[12pt]{article}

\usepackage{times}

\usepackage[a4paper,text={128mm,185mm}, centering]{geometry}

\usepackage{changepage}
\usepackage{titlesec}

\titleformat{\section}[hang]%
{\bfseries\large}{\thesection.}{1ex}{}%

\titleformat{\subsection}[hang]%
{\bfseries}{\thesubsection}{1ex}{}%

\usepackage{amsthm}
\usepackage{fancyhdr}
\theoremstyle{plain}
\newtheorem{theorem}{Theorem}%[section]

\newtheorem{corollary}[theorem]{Corollary}

\theoremstyle{definition}

\usepackage[matrix,arrow,curve,cmtip]{xy}
\usepackage{amssymb}
\usepackage{latexsym}
\usepackage{graphicx}
\usepackage{tikz}

\def\b0{{\bf 0}}
\def\b1{{\bf 1}}

\def\cC{{\cal C}}

\def\cB{{\cal B}}

\def\cR{{\cal R}}
\def\cC{{\cal C}}

\def\cS{{\cal S}}

\def\cQ{{\cal Q}}

\def\n{\noindent}

\def\eop{\hfill $\Box$ \smallskip}
\def\bop{\noindent {\bf Proof.} }

\title{\vskip 5pt  \bf  Canonical Sphere Bases for Simplicial and Cubical Complexes}
\author{\itshape\bfseries {Paul C. Kainen \thanks
{to appear in {\it Journal of Knot Theory \& its Ramifications}, accepted for publ. 6/25/22}
}}
\date{}

\begin{document}
\maketitle

%Nota bene:  the next two commands erase pages numbers (the right page numbers will be add by the editors on the pdf) 
%\cfoot{}
%\thispagestyle{empty}
%%
\vskip 25pt
% Abstract in French and in English, followed by Keywords and MSC:
\begin{adjustwidth}{0.5cm}{0.5cm}
{\small
%{\bf R\'esum\'e.} 

{\bf Abstract.} 
\n
Sphere-bases are constructed for the $\mathbb{Z}_2$ vector space formed by the $k$-dimensional subcomplexes, of $n$-simplex (or $n$-cube), for which every $(k{-}1)$-face is contained in a positive even number of $k$-cells;
addition is symmetric difference of the corresponding sets of $k$-cells.  The bases consist of the boundaries of an algorithmically-specified family of $k{+}1$-simplexes or $k{+}1$-cubes.  Geometric properties of these bases are investigated.
%The cardinality of the cube basis equals the rank of an integer homology group of a corresponding  {\it no-k-equal} space.
\\
{\bf Keywords.} Cellular tree, no $k$-equal space, robust and connected-sum  cycle bases for graphs, homology of Platonic skeleta, computational topology.\\
{\bf Mathematics Subject Classification (2010).} 57N99, 57M15, 55N31.
}
\end{adjustwidth}

% Here starts the text:.

\section{Introduction}

A cycle basis of a graph $G$ is a minimal family of connected 2-regular subgraphs such that each even-degree subgraph of $G$ is the mod-2 sum of a subset of the basis; ``mod-2'' means that an edge appears in the sum if and only if it appears an {\it odd} number of times in the sum.

Cycle bases of graphs have been constructed with special properties such as {\it fundamental} (constructable from a spanning tree) \cite{klemm-stadler}, {\it minimum}
(least total number of edges) \cite{horton}, and {\it robust} or {\it  connected sum} (having rearrangement properties similar to the ordering of facets in a shelling),  \cite{pck-2000, hk-2018}.  

We describe analogous bases for simplicial and cubical complexes with all or almost all of the above properties, suitably generalized to complexes. 

In addition, our bases are {\it canonical}: for $1 \leq k \leq n{-}1$, the basis for even $k$-dimensional subcomplexes of an $n$-dimensional simplex $\Delta_n$ or $n$-cube $\cQ_n$ is uniquely determined by 
the integer pair $(n,k)$ (together with a choice of basepoint for the simplex and a choice of coordinate-ordering for the cube) and could be precomputed and stored to speed up topological computations.  
Having canonical spheres as basis elements  might be useful for various aspects of persistent homology \cite{C-S-2007, dey-et-al, guerra-et-al,Z-C-2005}.

The general problem of realizing homology classes by spheres goes back to Steenrod \cite[Problem 24]{eilenberg-49}, who asked (as phrased by Eilenberg) {\it ``What algebraic conditions are necessary and sufficient for a homology class $z$ in  $H_{n}(K)$ with integer coefficients on a complex K to be spherical?''}   The problem has been extended to stable homology by Landweber \cite{landweber}.

With the convenient and natural properties of $\mathbb{Z}_2$ coefficients \cite[p. 7]{haus}, we get the $k$-th homology of cube and simplex $k$-skeleta spanned as vector spaces by our bases. The cross-polytope is not addressed here.

The simplex basis was known (though not its special property), while our cube basis seems to be new \cite{pck-amc}. 
The simplex basis for $k=1$ is robust \cite{pck-2000} and the cube basis has the connected sum property \cite{pck-amc}.  See \S5 below.
 
 Cardinality of each basis is the Betti number of the $k$-skeleton of the respective $n$-dimensional polytope (Section 2) and so expressible by means of the Euler-Poincar\'e equation.   Our constructions give alternative expressions which lead to binomial identities.   We also exhibit connections with topological combinatorics, cellular trees, and no-$k$-equal space (see \S 6).
 
For the $k$-skeleton of $\Delta_n$, our basis consists of the set of all boundaries of $k{+}1$-faces of $\Delta_n$ which contain a fixed vertex, so  the basis depends also on a base-point.   As boundary of a simplex is a sphere, this is a sphere basis.

In the $n$-cube, a sphere basis for the $k$-skeleton is obtained by a recursive construction which depends on the  sequential ordering of the coordinates.  

For instance, a basis $\cB = \{z_1, \ldots, z_{31}\}$ for the 2-skeleton of the 5-cube $[0,1]^5$ is given with respect to the ordering $(1,2,3,4,5)$ as follows: 
Let $$z_1 := \partial(Q(1,2,3)), \;\mbox{where}\;Q(1,2,3):=\{(x_1,\ldots,x_5): x_4 = 0 = x_5\}.$$\\
On each of the 6 square faces
of $Q(1,2,3)$, erect a 3-cube using the 4th coordinate, keeping the 5th coordinate fixed at zero. The boundaries of these six cubes are spheres $z_2, \ldots, z_7$, and $\{z_1,\ldots, z_7\}$ is a basis for the 2-skeleton of the 4-cube.  If we
now take each of the 24 square faces of the 4-cube and repeat this process using the 5th coordinate, we get 24 new spheres, and the resulting set of 31 spheres is a basis for the 2-skeleton of the 5-cube.  Order of the first three coordinates doesn't matter.  

Any even subgraph is an edge-disjoint union of cycles, so a cycle basis can be chosen for any graph such that any even subgraph {\it contains} the members of the basis of which it is the sum.  
But this is impossible for the  torus, an even 2-complex, which has no sphere subcomplex. However, the torus is a mod-$2$ sum of spheres in our basis. One can get a toroidal subcomplex in the 2-skeleton of $\cQ_4$, by summing five of the basis elements, excluding two ``opposite-side'' cube boundaries; see \cite[Fig. 6]{hk-monthly}. 
%is denoted$Q(a,b,\zeta)$, where $a$ and $b$  are two of the three coordinates and $\zeta$ is the fixed value of the third.  The boundaries of the cubes $$Q(a,b,\zeta,4) := Q(a,b,\zeta) \times [0,1].$$

%as illustrated for $\cQ_n^2$: $\cB(2,3)$ is the singleton $\{Q\}$ where $Q = \cQ_3^2$, while $\cB(2,4)$ has $7$ elements and consists of $Q$ together with the boundaries of the six $3$-cubes on the six sides of $Q$.  So $\cB(2,4)$ is a set of $7$ two-dimensional spheres. Every $2$-face of $\cQ_4$ is a face of some basis element.
 
%We encountered the two bases in studying a property of cycle bases of graphs.  

Section 2 below defines terms and lists results used. Section 3 is on the cube, and Section 4 is on the simplex. Section 5 generalizes cycle basis properties for $k>1$. The last section is a  discussion on coincidences in enumeration, binomial identities, and topics in geometric combinatorics.

\section{Complexes, cycle spaces, and homology}

We work with simplicial and cubical complexes, both of which are {\bf polytopal complexes} \cite[p. 127 ]{ziegler}.
For convenience, a few of the basic terms are defined; otherwise, see Ziegler \cite{ziegler} and Spanier \cite{Spanier}. We use $\mathbb{Z}_2$ for the 2-element vector space. For $d$ a positive integer,
$[d] := \{1,\ldots,d\}$;  $|A|$ denotes the cardinality of a set $A$.

The $n$-{\bf simplex} is the convex hull of the $n{+}1$ standard unit vectors 
\[
\Delta_n :={\rm conv}(\{e_1, \ldots, e_{n+1}\})
\]
and its $k$-{\bf faces} are the corresponding convex hulls of the various subsets 
$$\{e_j: j \in J \subseteq [n+1]\}, |J|=k{+}1.$$
The $n$-{\bf cube} is the product of $n$ copies of the unit interval,
$\cQ_n := [0,1]^n$.  Its $k$-{\bf faces} are the various subspaces of the form $A \times B$, where $A$ is the $k$-cube corresponding to some $k$-element subset of  $[n]$ and $B$ is a vertex in the $n{-}k$-cube determined by the complementary subset.

For 
$X := \Delta_n$ or $\cQ_n$,
let $X^k$ denote the $k$-{\bf skeleton} of $X$ (all faces of dimension $\leq k$).
A subcomplex $Y$ of $X$ has dimension $k$ if $k$ is the largest dimension of any face in $Y$ and we write $\dim Y := k$.
So $\dim X^k = k$.  We call $Y$ a {\bf k-complex} if $\dim Y = k$.  Note that every $k$-complex is a subcomplex of the $k$-skeleton of $X$ for $n$ sufficiently large. We say that $Y$ is a {\bf complex} if it is a  $k$-complex for some $k$.

For $\ell \geq 0$ and for $Y$ any $k$-complex contained in
$X$, we write $C_\ell(Y)$ for the $\mathbb{Z}_2$-vector space indexed by the set of all $\ell$-faces of $Y$, so $C_\ell(Y) = 0$ for $\ell > k$. There is a $\mathbb{Z}_2$-linear boundary operator $\partial_\ell: C_\ell(X) \to C_{\ell -1}(X)$ which maps an $\ell$-cell $c$ to the $(\ell{-}1)$-chain which is the formal sum of the $(\ell{-}1)$-faces of $c$, $\ell \geq 1$. Put $Z_\ell := Z_\ell(Y) := \ker \partial_\ell$ and $B_\ell := {\rm im}\, \partial_{\ell + 1} \subseteq Z_\ell$ and $H_\ell := H_\ell(Y)  := Z_\ell/B_\ell$. Let $b_\ell$ be the $\mathbb{Z}_2$-dimension of $H_\ell$ (Betti number).

A $k$-complex $Y$ is {\bf even} if each $(k{-}1)$-face in $Y$ is contained in a positive even number of $k$-faces in $Y$.  These are exactly the $k$-chains which map to zero by $\partial_k$; that is, $Z_k(Y)$ is the set of all even $k$-subcomplexes of $Y$.
Manifolds and pseudomanifolds \cite[pp. 148--150]{Spanier} are even complexes. 

For $X$ an $n$-complex and $k \leq n-1$, 
%There are no $(k{+}1)$-cells in $X^k$, for $1 \leq k < n$, so for $\mathbb{Z}_2$-homology
$$H_k(X^k)
%= Z_k(X^k)/B_k(X^k) := {\rm ker}(\partial_k)/{\rm im}(\partial_k) 
\cong Z_k(X^k) \cong Z_k(X).$$ 
Hence, $b_k(X^k)$ is the cardinality of a basis for the even $k$-complexes in $X$.

Each $\ell$-face of a complex $Y$ is an 
{\bf $\ell$-cell}, which is a closed topological ball of dimension $\ell$ (more specifically, a simplex or cube).  
Every (finite) complex has a unique topology coinduced by the topologies on its cells and is a compact Hausdorff space.  
For $T$ a topological space, a {\bf T-complex} is a complex which is homeomorphic to $T$.
%A $k$-complex is {\bf even} if each $(k{-}1)$-cell is in a positive even number of $k$-cells. 
Note that $X^k = \Delta^k_n$ or $ \cQ^k_n$ is even if and only if $n-k$ is odd as it is easy to verify that the number of $k$-cells containing any $k{-}1$-cell in both $n$-simplex and $n$-cube is $n-k+1$.

For $Y$ a $k$-complex, $Y^{(j)}$ denotes the set of all $j$-cells in $Y$, $j \geq -1$ with $Y^{(j)} = \emptyset$ for $j>k$; there is a unique ``empty'' cell in dimension $j = -1$. The elements in $Y^{(j)}$ are called {\bf vertices, edges, sides, {\rm or} facets} for $j = 0$, $1$, $k{-}1$, or $k$, resp.
A complex is {\bf pure} if every cell is a face of some facet.
%, and the {\bf $k$-skeleton} $Y^{k}$ of $Y$ is the subcomplex given by the cells of dimension $\leq k$. 

We use the {\it Mayer-Vietoris sequence}  \cite[p. 148]{haus}, which holds for an arbitrary field of coefficients (so for $\mathbb{Z}_2$).
For simplicial complexes $Y$, $Y'$, there is a long exact sequence in homology \cite[p. 187 ]{Spanier},
\begin{equation}
 \cdots \to H_{k}(Y \cap Y') \to H_{k}(Y) \oplus H_{k}(Y')\to H_{k}(Y \cup Y') \to H_{k-1}(Y \cap Y') \to \cdots.   
 \label{eq:mv}
\end{equation}
For $k \geq 2$, the middle terms in (\ref{eq:mv}) are isomorphic when $Y \cap Y'$ is contractible.
The same facts hold if $Y$ and $Y'$ are cubical cell complexes.

\section{Sphere bases for even subcomplexes of the cube}
  
Let $k \geq 2$.  For each $n \geq k+1$, we construct a basis $\cB(n,k)$ for 
%the even subcomplexes of the $k$-skeleton of the $n$-cube, 
$Z_k(\cQ^k_n)$, and then state its properties in a theorem. 

For $n=k {+} 1$, we have $\cB(n,k)= \{\partial \cQ_n\}$.  For $n \geq k{+}1$,
\begin{equation}
\cB(n{+}1,k) \setminus \cB(n,k) := \{\partial (s \times [0,1]): s \in \cQ_n^{(k)}\} 
\label{eq:incr}
\end{equation}
It is the linear ordering of coordinates (after the first $k{+}1$) that determines which $k{+}1$-cube boundaries are in the basis.
The recursively generated set $\cB(n,k)$ of spherical $k$-dimensional subcomplexes is independent over $\mathbb{Z}_2$. Indeed,  each additional isomorphic copy of $\partial \cQ_{k+1}$, after the first, contains a $k$-face (opposite to the $k$-face to which the copy of $\cQ_{k+1}$ is attached) which is in no other sphere, so a linear combination is zero if and only if all coefficients of combination are also zero.  We now prove $\cB(n,k)$ is a basis.

\begin{theorem}
Let $k \geq 1$.  Then for every $n \geq k{+}1$,  (i) $\cB(n,k)$ is a basis for $Z_k(\cQ_n^k)$ and (ii) every $k$-cube in $\cQ_{n+1}$ is a face of some member of $\cB(n,k)$.
\label{th:cubeBasis}
\end{theorem}
\begin{proof}
The result was proved for $k=1$ in \cite{pck-amc}; we now fix $k \geq 2$ and proceed by induction on $n$; trivially, (i) and (ii) hold if $n=k+1$.  
By the inductive hypothesis, for $n \geq k+1$, 
$|\cB(n,k)| = b_k(\cQ_n^k)$ and 
$$U(\cB(n,k)) := \{c: c \in  S\;\mbox{and} \; S \in \cB(n,k)\} = \cQ_n^{(k)}.$$
Then (ii) holds for $n+1$ as
$$U(\cB(n+1,k) = \cQ_n^{(k)} \cup \{c: c \in \partial (s \times [0,1]), s \in \cQ_n^{(k)}\} = \cQ_{n+1}^{(k)}.$$
Hence, $\cQ_{n+1}^k$ is the union of $\cQ_n^k$ and the 
$k$-spheres in $\cB(n+1,k) \setminus \cB(n,k)$, which we call the {\it attached} spheres.
In fact, for each pair $(n,k)$ statement (ii) follows from (i) but we actually need (ii) in order to prove (i).
 
For any $\ell \geq 2$, the $\ell{+}1$-cube is the (Cartesian) product of an $\ell$-cube and the unit interval.  The {\it back} $\ell$-cube corresponds to $0$ in the last coordinate, while  {\it front} $\ell$-cube corresponds to $1$.  The {\it side} $\ell$-cubes are the products of the unit interval with the $\ell{-}1$-subcubes of the original $\ell$-cube.
% and the unit interval.  
Order the sequence of ${{n}\choose{k}}2^{n-k}$ attached spheres arbitrarily, $S_1, S_2, \ldots, S_t$; then $S_i$ 
intersects the union of $\cQ_n^k$ and $\bigcup_{r=1}^{i-1}S_r$ in its back $k$-face together with a subset of the side cubes.   
Hence, the intersection of each sphere with the union of $\cQ_n^k$ and all previously attached spheres  
is contractible. 

For each $i$, $1 \leq i \leq t = {{n}\choose{k}}2^{n-k}$, applying the Mayer-Vietoris isomorphism (\ref{eq:mv}) with $Y$ equal to $\cQ_n^k \cup \bigcup_{r=1}^{i-1} S_r$ and $Y' = S_i$ increments the dimension of the $k$-th homology of the union by $1$
for each attachment of a sphere .  Hence, $b_k(\cQ_{n+1}^k) =  b_k(\cQ_{n}^k) + {{n}\choose{k}}2^{n-k}$, so
using (i) for $\cQ_n^k$ and (\ref{eq:incr}),
$$b_k(\cQ_{n+1}^k) 
= |\cB(n,k)| + {{n}\choose{k}}2^{n-k} = |\cB(n+1,k)|.$$ 
Thus, for every $n$, $\cB(n,k)$ is a maximal independent set and so is a basis. 
\end{proof}

\section{Sphere bases for even subcomplexes of the simplex}
%{The cycle space of $\Delta_d$ has a basis of spheres}

%The situation for the simplex is simpler than for the cube.

For $1 \leq k \leq n-1$, 
%let $m'(n,k)$ denote the alternating sum of the number of $k$-simplices contained in  $\Delta_n$. B
by Euler-Poincar\'e,  $b_k(\Delta_n^k) = m'(n,k)$, where
\begin{equation}
m'(n,k) := \sum_{j=-1}^k (-1)^{k-j}{{n+1}\choose{j+1}} = {{n+1}\choose{k+1}} -{{n+1}\choose{k}} \pm \cdots +(-1)^{k+1}.
\label{eq:betti-simplex}
\end{equation}

%$m'(2,4) = -1 + 5 - 10 + 10 ={{4}\choose{3}}$.
%while $m'(3,7) = 1 - 8 + 28 - 56 + 70 = 35 = {{7}\choose{4}}$.

%For $1 \leq k \leq n-1$, l

Let $\cC'(n,k)$ be the family of all $k{+}1$-simplices 
$\Delta'$ 
contained in $\Delta_n$ and containing vertex  $1$ of $\Delta_n$. 
%But these are just the joins of the vertex $d+1$ with the $k$-simplexes in $\Delta_{d-1}$. 
Put
$\cB'(n,k) := \{\partial c: c \in \cC'(n,k)\}$.  By definition
\begin{equation}
|\cB'(n,k)| = {{n}\choose{k+1}}.
\label{eq:kbasis-simplex}
\end{equation}

\begin{theorem}
Let $1 \leq k \leq n{-}1$. Then $\cB'(n,k)$ is a sphere basis for $\mathbb{Z}_{k}(\Delta_n)$.
\label{th:sx-k-basis}
\end{theorem}
\bop 
Write
${{n}\choose{k+1}} - m'(n,k)$ as a sum of $k{+}3$ terms and use the binomial recursion $k{+}1$ times to get $(-1)^{k+2} + (-1)^{k+1} =0$. Hence, $\cB'(n,k)$ has the cardinality of a basis. Because it is independent, it must {\it be} a basis.
%since each element has a $k$-simplex in none of the others.
%$${{n}\choose{k+1}} - {{n+1}\choose{k+1}} + {{n+1}\choose{k}} - {{n+1}\choose{k-1}} + - \cdots + (-1)^{k+1} =$$
%$$= -{{n}\choose{k}} + {{n+1}\choose{k}} - {{n+1}\choose{k-1}} \pm \cdots + (-1)^{k+1} = \cdots 
% = (-1)^{k+2} + (-1)^{k+1} 
%=0.$$
\eop
%$$= {{n}\choose{k-1}} - {{n+1}\choose{k-1}}+ - \cdots + (-1)^{k+1} =$$ 
%$$= \cdots = (-1)^{k+2}{{n}\choose{0}} + (-1)^{k+1} = (-1)^{k21} + (-1)^{k+1} = 0$$
%using the basic binomial recursion repeatedly. 

As a consequence, $b_k(\Delta_n^k) = {{n}\choose{k+1}}$; cf. \cite[\S 4.1]{dkm-2016}, \cite{stack}. 
%Indeed, this is well-known that the Betti number of the $k$-skeleton of an $n$-simplex is given by (\ref{eq:kbasis-simplex})  (e.g., \cite{stack}). 
The simplicial complex generated by the members of $\cB'(n,k)$ is a {\bf shifted} complex \cite{ssad}, so $b_k(\Delta_n^k)$ is the number of facets not containing $1$, agreeing with our result. 
A pure $k$-complex is {\bf shellable} if it is possible to order its facets so that each facet intersects the union of the previous facets in a pure $k{-}1$-complex.  
Sphere bases exist for shellable complexes \cite{bjorner-1990}, the boundary complex of a polytope is shellable \cite{bm72}, and the $k$-skeleton of a shellable complex is shellable \cite{bw96}.  The basis $\cB'(n,k)$  arises also in  
matroid theory 
\cite{bjorner-1990}.
% as the simplex is a {\it matroid independence complex}.

\section{Some properties of the bases}

We first note the consequence of our two theorems.
\begin{corollary} Let $Y$ be a pure, even $k$-subcomplex of $X = \Delta_n$ or $\cQ_n$. Then $Y$ is the mod-$2$ sum of a family of at most $b_k(X^k)$ distinct $k$-dimensional spheres in the respective bases $\cB'(n,k)$ or $\cB(n,k)$. \label{co:simplex-cx}\end{corollary}

The property of {\it minimality} for a cycle basis certainly applies to both of our bases since each uses only minimum-size elements in $Z_k(X)$, where $X$ is $n$-simplex or $n$-cube.  Any spanning tree in a connected graph determines a unique cycle basis, where each cycle is produced by one of the non-tree edges of the graph.  Such cycle bases are called {\it fundamental}; interpreting this property for $k>1$ depends on generalizing ``spanning tree''.

Both bases $\cB'(n,k)$ and $\cB(n,k)$, when $k=1$, have recursive properties not possessed by all cycle bases \cite{hk-2018}.  It is possible that some similar properties hold for $k > 1$ since the key condition in the recursion, that two cycles being added meet in a common nontrivial path, generalizes -- the connected sum of two $k$-spheres remains a $k$-sphere. But $\cB, \cB'$ aren't panaceas \cite{hammack}.

A cycle basis $\cR$ for a graph is {\bf robust} if for each cycle $z$ in the graph, it is possible to find a (finite) sequence of members of $\cR$ such that $z$ is their mod-2 sum and each summand intersects the sum of the previous terms in a nontrivial path \cite{pck-2000}.  Hence, each partial sum is also a cycle. 

The cycle basis $\cB'(n,1)$ is a robust basis for $Z_1(K_{n+1})$ \cite[Prop. 1]{pck-2000}.  
Extending the notion  to $k>1$, we ask: Is $\cB'(n,k)$ a robust basis for $Z_k(\Delta_n)$?

The cycle basis $\cB(n,1)$ for the hypercube graph $Q_n$ has a weaker recursive property, called {\bf connected sum} \cite{pck-amc, hk-2018}:  Any cycle in the graph can be constructed by iterating the procedure described for robust bases.  
The first iteration constructs a family of cycles from the original basis.  In the later iterations, one uses both the basis and the cycles constructed in the previous iterations, until, after finitely many iterations, all cycles are formed.   
Extending this to $k>1$, we ask if $\cB(n,k)$ is a connected sum basis for $Z_k(\cQ_n)$?

A property of cycles is called {\bf cooperative} \cite{hk-2018} if the mod-2 sum of two cycles, which intersect in a nontrivial path, has the property whenever both the summands have it.
%A virtue of connected sum and robust bases is that 
{\it Any cooperative property which holds for all cycles in a connected-sum basis must hold for all cycles in the graph}.  For example, $Q_5$ has over $51$ billion cycles (sequence A085408 in the OEIS \cite{oeis}) but has only  $49$ elements in the cycle basis $\cB(5,1)$ which is a connected sum basis.

These bases provide spherical primitives which could be added together sequentially to build up any spherical object.  Cubical complexes seem more natural for applications, especially in graphics, but 
consider a $k$-dimensional  simplicial complex $Y$ %\cite[p. 148]{Spanier}, 
and let $\Delta_n$ be a simplex which contains $Y$ and has at least one additional vertex $1$.  Let 
 $1*\sigma$ denote the $(k{+}1)$-simplex which is the topological join of $1$ and a $k$-simplex $\sigma$ in $Y$. Then we define
\begin{equation}
M(Y) := M := \sum_{\sigma \in Y^{(k)}}\partial(1 \,*\,\sigma)
\end{equation}
so $M$ is the mod-2 sum of these canonically given members of $\cB'(n,k)$.

This may be an advantage for simplicial, rather than cubical, models.

An entirely different approach to the geometry of $k$-skeleta of Platonic polytopes was taken in \cite{hk-2018-platonic, hk-2020, hkgeomb}, where we were able to {\it decompose} the even $k$-skeleta of cube and simplex into {\it facet-disjoint} spheres in nearly all cases for cubes (and always for simplexes).  

\section{Discussion}

Fix $k \geq 1$ and let $n \geq k{+}1$.  Our construction shows
%the integer sequences $s(n,k) := |\cB(n,k)|$, for $n \geq k+1$, by the construction of the basis and equation (\ref{eq:incr}), can be expressed
\begin{equation}
s(n,k) := |\cB(n,k)| = \sum_{j=k}^{n-1}{{j}\choose{k}} 2^{j-k}.
\label{eq:s}
\end{equation}
 But $m(n,k) := b_k(\cQ_n^k)$ can  be computed via 
the Euler-Poincar\'e formula
\cite[p. 146]{Hatcher}, \cite[p. 25]{haus}, to get a rather different-looking expression.
% Euler-Poincar\'e  
\begin{equation}
 m(n,k) = (-1)^{1+k}+\sum_{j=0}^k (-1)^{k-j}{{n}\choose{j}} 2^{n-j}.
\label{eq:betti-cube}
\end{equation}

Theorem \ref{th:cubeBasis} implies that $s(n,k) = m(n,k)$.
Conversely, a computer-algebra proof that $s(n,k) = m(n,k)$ implies Theorem \ref{th:cubeBasis}.
% as the spheres in $\cB(n,k)$ are  independent. 
Indeed, both $m$ and $s$ satisfy the recursion $ T(n,k) = 2\,T(n-1,k) + T(n-1,k-1)$
for $1 \leq k \leq n-1$
%: $\;\;$
%\begin{equation} T(n,k) = 2\,T(n-1,k) + T(n-1,k-1).    \label{eq:recur} \end{equation}
as can be verified using Mathematica \cite{matha}, which applies a Gauss hypergeometric contiguous identity - see the Digital Library of Mathematical Functions \cite[15.5.15]{dlmf}.

The sequences given by $s(n,k)$ for $n=k{+}1, k{+}2, \ldots$ with  $k = 1, 2, 3, \ldots$ \\
$1, 5, 17, 49, 129, 321, 769, 1793, \ldots$\\
$1, 7, 31, 111, 351, 1023 , 2815, 7423, \ldots$\\
$1, 9, 49, 209, 769, 2561, 7937, 23297, \ldots $\\
constitute a  triangular array; see the OEIS \cite[A119258]{oeis} and   \cite{sw-2010}.

For $1 \leq k \leq n{-}1$, the sequences are the ranks of certain (integer) cohomology and homology groups \cite{bly92, bw-95, bary-2017} associated with the {\bf no k-equal} space $M^{\mathbb{R}}_{n,k}$,  which is the complement in $\mathbb{R}^n$ of 
$$V^{\mathbb{R}}_{n,k} := \{x \in \mathbb{R}^n : \exists  J \subseteq [n]\;s.t. \;|J|=k\;\mbox{and}\,\;x_j\;\mbox{{\rm is constant $\forall j \in J$}}\}.$$
In particular, Bjorner and Welker found \cite{bw-95} that for $n \geq k$,
\begin{equation}
{\rm bw}(n,k):={\rm rank}\Big(H^{k-2}(M^{\mathbb{R}}_{n,k}, \mathbb{Z})\Big) = \sum_{i=k}^n {{n}\choose{i}}\,{{i-1}\choose{k-1}}, \,k \geq 3.
\label{eq:bw}
\end{equation}
Evaluation shows that ${\rm bw}(n,3) = s(n,2)$ is sequence A055580 from \cite{oeis}, while ${\rm bw}(n,4) = s(n,3)$ is A027608, and the next sequence is A211386.

There are a number of distinct formulas that all give these same sequences.  See  \cite{half-cube, green-harper, sw-2010} which prove the formulas are equal algebraically. 

%Green \cite{half-cube} also found a homological connection with 
Green \cite[Thm. 4.1.2]{half-cube}  shows that, for a CW-complex $C_{n,k}$ contained in the  {\it half-cube} \cite[\S 8.6] {coxeter}, $b_{k-1}(C_{n,k})$ can be expressed by the alternating sum
%by deleting the interiors of all ``half-cube shaped'' $\ell$-cells for $\ell \geq k$.   He then shows that
\begin{equation}
b_{k-1}(C_{n,k}) = {\rm gr}(n,k) := \sum_{i=k}^n (-1)^{k+i}\, 2^{n-i}{{n}\choose{i}}
\label{eq:gr}
\end{equation}
and Green also shows \cite[Cor. 4.1.6]{half-cube}
\begin{equation}
{\rm gr}(n,k) = {\rm bw}(n,k).
\end{equation}

Tree-like structures in higher dimensions have been studied for decades; e.g., Pippert \&  Beineke \cite{pb-1969}, Dewdney \cite{dewdney}, Bolker \cite{bolker}, Kalai \cite{kalai}, and Lyons \cite{lyons}.
Duvall, Klivans, and Martin \cite{dkm-2016} define a  {\it $k$-dimensional cellular spanning tree} of a $k$-complex to be a certain subset of the facets containing the $n{-}1$-skeleton and satisfying some conditions on integer homology.  

Baryshnikov, Klivans, and Kosar \cite[Thm. 1.1]{bary-2017} show that the number of facets in a $k{+}1$-cellular spanning tree $T$ for $\cQ_n^{k+1}$ is given by $b_{k-1}(M^{\mathbb{R}}_{n,k-1})$.  Further, if $T$ is a spanning tree of the $k{+}1$-skeleton of {\it any convex polytope} $P$ in $\mathbb{R}^n$, they show that $T$ has $|T|=b_k(P^k)$ facets
%satisfies \begin{equation}|T^{(k+1)}| = b_{k-1}(P^{k})\end{equation} $b_{k-1}(P^{k-1})$.
\cite[Prop. 5.2]{bary-2017}.

Let $\cC(n,k)$ be the family of $k{+}1$-cubes whose boundaries are $\cB(n,k)$.
We conjecture that $\cC(n,k)$  is a $(k{+}1)$-cellular spanning tree for $\cQ_n^{k+1}$ and that $\cC'(n,k)$  is a $(k{+}1)$-cellular spanning tree for $\Delta_n^{k+1}$.

%\subsubsection*{Acknowledgements}

{\bf Acknowledgements:}
I thank colleagues in  Slovenia, (Drago Bokal) and Greece (Stratos Prassidis and Sofia Lambropoulou)
for a stimulating research environment during the preparation of the first draft of this paper in the summer of 2018.
Thanks also to the referees for helpful comments and for pointing out the matroid connection and to Michael Somos for observing the role of the Gauss identity.

% Cahiers wants the author's address at the end of the paper:

\vspace{5mm}
\noindent
Paul C. Kainen,
%Department of Mathematics and Statistics
Georgetown University\\
%37th and O Streets, N.W.\\
%Washington, D.C. 20057 USA\\
kainen@georgetown.edu

\end{document}